\documentclass[reqno]{amsart}
\usepackage{graphicx}
\usepackage{pictex}
\usepackage{mathrsfs}
\usepackage{hyperref}
\begin{document}

\newtheorem{thm}{Theorem}[section]
\newtheorem{lem}[thm]{Lemma}
\newtheorem{cor}[thm]{Corollary}
\newtheorem{add}[thm]{Addendum}
\newtheorem{prop}[thm]{Proposition}
\newtheorem{conj}[thm]{Conjecture}
\theoremstyle{definition}
\newtheorem{defn}[thm]{Definition}

\theoremstyle{remark}
\newtheorem{rmk}[thm]{Remark}


%
%
%
%

\newcommand{\OmegaH}{\Omega/\langle H \rangle}
\newcommand{\hatOmegaHstar}{\hat \Omega/\langle H_{\ast}\rangle}
\newcommand{\SurfG}{\Sigma_g}
\newcommand{\TriangG}{T_g}
\newcommand{\TriangGOne}{T_{g,1}}
\newcommand{\ProjG}{\mathcal{P}_g}
\newcommand{\TeichG}{\mathcal{T}_g}
\newcommand{\CirclePackGTau}{\mathsf{CPS}_{g,\tau}}
\newcommand{\CrossRatio}{{\bf c}}
\newcommand{\CrossRatioGTau}{\mathcal{C}_{g,\tau}}
\newcommand{\CrossRatioOneTau}{\mathcal{C}_{1,\tau}}
\newcommand{\DeformGTau}{\mathcal{C}_{g,\tau}}
\newcommand{\Forget}{\mathit{forg}}
\newcommand{\Uniform}{\mathit{u}}
\newcommand{\Section}{\mathit{sect}}
\newcommand{\SLTwoC}{\mathrm{SL}(2,{\mathbb C})}
\newcommand{\SLTwoR}{\mathrm{SL}(2,{\mathbb R})}
\newcommand{\SUTwo}{\mathrm{SU}(2)}
\newcommand{\PSLTwoC}{\mathrm{PSL}(2,{\mathbb C})}
\newcommand{\GLTwoZ}{\mathrm{GL}(2,{\mathbb Z})}
\newcommand{\GLTwoC}{\mathrm{GL}(2,{\mathbb C})}
\newcommand{\PSLTwoR}{\mathrm{PSL}(2,{\mathbb R})}
\newcommand{\PGLTwoR}{\mathrm{PGL}(2,{\mathbb R})}
\newcommand{\PSLTwoZ}{\mathrm{PSL}(2,{\mathbb Z})}
\newcommand{\SLTwoZ}{\mathrm{SL}(2,{\mathbb Z})}
\newcommand{\PGLTwoZ}{\mathrm{PGL}(2,{\mathbb Z})}
\newcommand{\nnn}{\noindent}
\newcommand{\MCG}{{\mathcal {MCG}}}
\newcommand{\MMap}{{\bf \Phi}_{\mu}}
\newcommand{\HH}{{\mathbb H}^2}
\newcommand{\TT}{{\mathbb T}}
\newcommand{\X}{{\mathcal  X}}
\newcommand{\CC}{{\mathbb C}}
\newcommand{\RR}{{\mathbb R}}
\newcommand{\Q}{{\mathbb Q}}
\newcommand{\ZZ}{{\mathbb Z}}
\newcommand{\PL}{{\mathscr {PL}}}
\newcommand{\GP}{{\mathcal {GP}}}
\newcommand{\GT}{{\mathcal {GT}}}
\newcommand{\GQ}{{\mathcal {GQ}}}
\newcommand{\EE}{{{\mathcal E}(\rho)}}
\newcommand{\HHH}{{\mathbb H}^3}

\title{The ${ \SLTwoC}$ character variety of the one-holed torus }
\author{Ser Peow Tan}
\address{Department of Mathematics \\ National University of Singapore \\
2 Science Drive 2 \\ Singapore 117543} \email{mattansp@nus.edu.sg}
\thanks{The authors are partially supported by the National University
of Singapore academic research grant R-146-000-056-112. The third
author is also partially supported by the National Key Basic
Research Fund (China) G1999075104.}

\author{Yan Loi Wong}
\address{Department of Mathematics \\ National University of Singapore \\
2 Science Drive 2 \\ Singapore 117543} \email{ matwyl@nus.edu.sg}

\author{Ying Zhang}
\address{Department of Mathematics \\ National University of Singapore \\
2 Science Drive 2 \\ Singapore 117543; current address: Department
of Mathematics \\Yangzhou University \\Yangzhou 225002
\\P. R. China} \email{yingzhang@alumni.nus.edu.sg
}

 \commby{}
\date{}

%

%
%

\begin{abstract}
In this note we announce several results concerning the $\SLTwoC$
character variety ${\mathcal X}$ of the one-holed torus. We give a
description of the largest open subset ${\mathcal X}_{BQ}$ of
${\mathcal X}$  on which the mapping class group $\Gamma$ acts
properly discontinuously, in terms of two very simple conditions,
and show that a series identity generalizing McShane's identity
for the punctured torus holds for all characters in this subset.
We also give variations of the McShane-Bowditch identities to
characters fixed by an Anosov element of $\Gamma$ with
applications to closed hyperbolic three manifolds. Finally we give
a definition of end invariants for $\SLTwoC$ characters and give a
partial classification of the set of end invariants of a character
in ${\mathcal X}$.
\end{abstract}

\maketitle

\vskip 20pt
\section{{\bf Introduction}}\label{s:intro}
\vskip 20pt

Let $T$ be the one-holed torus, $\pi$ its fundamental group, and
$\Gamma:=\pi_0({\rm Homeo}(T))$ the mapping class group of $T$. In
this note we announce several results concerning the $\SLTwoC$
character variety ${\mathcal X}$ of $T$.
Here is a brief description of our results. We first give a
characterization of the largest open subset ${\mathcal X}_{BQ}$ of
${\mathcal X}$ on which the mapping class group $\Gamma$ acts
properly discontinuously, in terms of two very simple conditions,
called the Bowditch Q-conditions. This generalizes results of
Bowditch \cite{bowditch1998plms}, who gives a similar description
for the ``type-preserving'' characters, and also of Goldman
\cite{goldmanGT2003}, see also \cite{goldmanAM1997}, who studied
the dynamics of the action of $\Gamma$ on the \emph{real} ${\rm
SL}(2)$ characters, and gave a (geometric) description of the set
for the real characters. Note that for $[\rho]\in {\mathcal X}$,
it is possible to verify if $[\rho]$ satisfies these conditions
algorithmically.

We next show that a series identity generalizing McShane's
remarkable identity for the punctured torus holds for all
characters in this subset (the original identity can be regarded
as the formal derivative of the general identity evaluated at the
appropriate parameter value). This generalizes results of McShane
\cite {mcshane1991thesis}, \cite{mcshane1998im}, Bowditch
\cite{bowditch1996blms}, \cite{bowditch1998plms}, Mirzakhani
\cite{mirzakhani2004preprint} and the authors
\cite{tan-wong-zhang2004cone-surfaces},
\cite{tan-wong-zhang2004schottky} for two generator subgroups of
$\SLTwoC$. We also give necessary and sufficient conditions for
this identity to hold for characters $[\rho]\in {\mathcal X}$,
thereby giving a complete answer to the question of when the
identity holds for two generator subgroups of $\SLTwoC$.

The next set of results are for variations of the McShane-Bowditch
identities to characters fixed by an Anosov element of $\Gamma$,
and which satisfy a relative version of the Bowditch Q-conditions.
These have applications to closed hyperbolic three manifolds, and
generalize the results of Bowditch in \cite{bowditch1997t}.

Finally we give a definition of end invariants for $\SLTwoC$
characters, inspired by Bowditch's definition in
\cite{bowditch1998plms},  and give a partial classification of the
set $\EE$ of end invariants of a character  $[\rho] \in{\mathcal
X}$. When non-empty, this set gives information about the extent
to which the action of $\Gamma$ on $[\rho]$ is not proper. In
particular we give classification results for real characters,
imaginary characters and discrete characters.

\vskip 5pt

\noindent {\it Acknowledgements.} We would like to thank Bill
Goldman, Caroline Series, Makoto Sakuma  and Greg McShane for
their encouragement, helpful conversations, correspondence and
comments.

\vskip 20pt
\section{{\bf Preliminaries and definitions}}\label{s:preliminaries}
\vskip 20pt

\subsection{Basic Definitions}\label{ss:definitions}
Let $T$ be the one-holed torus and $\pi$ its fundamental group
which is freely generated by two elements $X,Y$ corresponding to
simple closed curves on $T$ with geometric intersection number
one. The $\SLTwoC$ character variety ${\mathcal X}:={\rm Hom}(\pi,
\SLTwoC)/\!/\SLTwoC$ of $T$ is the set of equivalence classes of
representations $\rho:\pi \mapsto \SLTwoC$, where the equivalence
classes are obtained by taking the closure of the orbits under
conjugation by $\SLTwoC$. The character variety stratifies into
relative character varieties: for $\kappa \in \CC$, the $\kappa$
relative character variety ${\mathcal X}_{\kappa}$ is the set of
equivalence classes $[\rho] $ such that
$${\rm tr}\,\rho(XYX^{-1}Y^{-1})=\kappa$$ for one (and hence any)
pair of generators $X,Y$ of $\pi$. By classical results of Fricke,
we have the following identifications:

$${\mathcal X} \cong \CC^3,$$
$${\mathcal X}_{\kappa} \cong \{(x,y,z)\in \CC^3 ~|~ x^2+y^2+z^2-xyz-2=\kappa\},$$
the identification is given by
$$\iota:[\rho] \mapsto (x,y,z):=({\rm
tr}\rho(X), {\rm tr}\rho(Y), {\rm tr}\rho(XY)),$$ where $X, Y$ is
a fixed pair of generators of $\pi$. The topology on ${\mathcal
X}$ and ${\mathcal X}_{\kappa}$ will be that induced by the above
identifications.

A character is \emph{real} if $\iota([\rho]) \in \RR^3$,
\emph{imaginary} if two of the entries of $\iota([\rho])$ are
purely imaginary and the third real, and dihedral if two of the
entries of $\iota([\rho])$ are zero (so that the third entry is
$\pm \sqrt{\kappa+2}$).

The outer automorphism group of $\pi$, ${\rm Out}(\pi):={\rm
Aut}(\pi)/{\rm Inn}(\pi)\cong {\rm GL}(2, \mathbb Z) $ is
isomorphic to the mapping class group $\Gamma:=\pi_0({\rm
Homeo}(T))$ of $T$ and acts on ${\mathcal X}$ preserving the trace
of the commutator of a pair of generators, hence it also acts  on
${\mathcal X}_{\kappa}$, the action is given by
$$\phi([\rho])=[\rho \circ \phi^{-1}],$$ where $\phi \in {\rm
Out}(\pi)$ and $[\rho] \in {\mathcal X}$ or ${\mathcal
X}_{\kappa}$ respectively. It is often convenient to consider only
the subgroup ${\rm Out}(\pi)^+$ of ``orientation-preserving''
automorphisms, corresponding to the orientation-preserving
homeomorphisms $\Gamma^+$ of $T$, which is isomorphic to
$\SLTwoZ$. The action of ${\rm Out}(\pi)^+$ (respectively, ${\rm
Out}(\pi)$) on ${\mathcal X}$ and ${\mathcal X}_{\kappa}$ is not
effective, the kernel is $\{\pm I\}$, generated by the elliptic
involution of $T$ so that the effective action is by $\PSLTwoZ$
(respectively, $\PGLTwoZ$).

\vskip 5pt \subsection{Simple curves, pants graph of
$T$}\label{ss:simplecurves}
 Let ${\mathscr C}$ be the set of free homotopy classes
of non-trivial, non-peripheral simple closed curves on $T$, the
elements of ${\mathscr C}$ correspond to certain elements of
$\pi$, up to conjugation and inverse. Let ${\mathscr C}(T)$ be the
``pants graph'' of $T$, the vertices of ${\mathscr C}(T)$ are the
elements of ${\mathscr C}$, where two vertices are joined by an
edge
 if and only if the corresponding curves on $T$
have geometric intersection number one. $\Gamma$ and ${\rm
Out}(\pi)$ acts naturally on ${\mathscr C}$ (respectively
${\mathscr C}(T)$). We can realize ${\mathscr C}(T)$ as the Farey
graph/triangulation of the upper half plane $\HH$ so that
${\mathscr C}$ is identified with $\hat \Q$, the action of
$\Gamma$ is realized by the action of $\PGLTwoZ$ on the Farey
graph. The projective lamination space $\PL$ of $T$ is then
identified with $\hat \RR$ and contains ${\mathscr C}$ as the
(dense) subset of rational points.

\vskip 5pt

\subsection{Bowditch Q-conditions (BQ-conditions)}\label{ss:BQconditions}
For $[\rho] \in {\mathcal X}$ and $X \in {\mathscr C}$, ${\rm
tr}\, \rho(X)$ is well-defined. We define the \emph{Bowditch
space} as the subset ${\mathcal X}_{BQ}$ of ${\mathcal X}$
consisting of characters $[\rho]$ satisfying the following
conditions (the \emph{Bowditch Q-conditions}):
\begin{enumerate}
    \item ${\rm tr}\,\rho(X) \not \in [-2,2]$ for all $X \in \mathscr C$;
    \item $|{\rm tr}\,\rho(X)| \le 2$ for only finitely many (possibly no)
    $X \in \mathscr C$.
\end{enumerate}

\vskip 5pt

\section{Statement of results}\label{s:Statements}

\subsection{Quasi-convexity}\label{ss:quasiconvexity}
For fixed $[\rho] \in {\mathcal X}$ and $K>0$, define ${\mathscr
C}_K(T)$ to be the subgraph of ${\mathscr C}(T)$ spanned by the
set of $X \in \mathscr C$ for which   $|{\rm tr} \rho(X)| \le K$.
Then we have the following simple but fundamental result of
Bowditch which plays a key role in the proofs of all subsequent
results. The result was stated for type-preserving characters in
\cite{bowditch1998plms} but the proof works for all characters.

\begin{thm}\label{thm:quasiconvexity}{\rm (}Bowditch
\cite{bowditch1998plms}{\rm )} For any $[\rho] \in {\mathcal X}$,
${\mathscr C}_K(T)$ is connected for all $K \ge 2$.
\end{thm}

\subsection{Action of $\Gamma$ and generalizations of McShane's
identity}\label{action of gamma}

\begin{thm} \label{thm:TWZ}{\rm (}Theorems 2.2, 2.3 and Proposition 2.4 of
\cite{tan-wong-zhang2004gMm}{\rm )}
\begin{enumerate}
    \item [(a)] ${\mathcal X}_{BQ}$ is open in ${\mathcal X}$.
    \item [(b)] $\Gamma$  acts properly discontinuously on ${\mathcal
    X}_{BQ}$. Furthermore, ${\mathcal X}_{BQ}$ is the largest open
    subset of ${\mathcal X}$ for which this holds.
    \item [(c)] For $[\rho] \in {\mathcal X}_{BQ} \cap {\mathcal X}_{\kappa}$,
    \begin{eqnarray}\label{eqn:TWZ}
\sum_{X \in {\mathscr
C}}\log\frac{e^{\nu}+e^{l(\rho(X))}}{e^{-\nu}+e^{l(\rho(X))}} =\nu
\mod 2 \pi i,
\end{eqnarray}
and the sum converges absolutely; where
$\nu=\cosh^{-1}(-\kappa/2)$, and the complex length $l(\rho(X))$
is given by the formula
\begin{equation}\label{eqn:complexlength}
l(\rho(X))=2\cosh ^{-1}(\frac{tr~ \rho(X)}{2}),
\end{equation}
 where we adopt the convention that the function $\cosh^{-1}$ has
images with non-negative real part and imaginary part in $(-\pi,
\pi]$.
\end{enumerate}

\end{thm}

\begin{rmk}\label{rmk:1}
In the case when $\kappa=-2$, $\nu=0$ and all the terms of
(\ref{eqn:TWZ}) are identically zero. However, if we take the
first order infinitesimal, or the formal derivative of
(\ref{eqn:TWZ}) and evaluate at $\nu=0$, we get
\begin{equation}\label{eqn:McShane}
\sum_{X \in {\mathscr C}}\frac{1}{1+e^{l(\rho(X))}} =\frac 12,
\end{equation}
which is McShane's original identity in \cite{mcshane1991thesis}
for real type-preserving characters, and also Bowditch's
generalization in \cite{bowditch1996blms} and
\cite{bowditch1998plms} for general type-preserving characters.
Note also that when $\kappa=2$, which corresponds to the reducible
characters,  the Bowditch Q-conditions are never satisfied, see
\cite{tan-wong-zhang2004endinvariants}. Parts (a) and (b) of the
above were originally stated in \cite{tan-wong-zhang2004gMm} in
terms of the relative character varieties.
\end{rmk}

\subsection{Necessary and sufficient conditions}\label{ss:necandsuff}
 Replacing condition (1) of the BQ-conditions by ($1'$) ${\rm
tr}\,\rho(X) \not \in (-2,2)$ for all $X \in \mathscr C$, we get
the \emph{extended Bowditch space} $\hat {\mathcal X}_{BQ}$, and
we have  the following result:
\begin{thm} \label{thm:necsuff} {\rm (}Theorem 1.6 of
\cite{tan-wong-zhang2004necsuf}{\rm )} For $[\rho] \in {\mathcal
X}$, the identity (\ref{eqn:TWZ}) of Theorem \ref{thm:TWZ} (c)
holds (with absolute convergence of the sum) if and only if
$[\rho]$ lies in the extended Bowditch space $\hat {\mathcal
X}_{BQ}$.
\end{thm}

The above result gives a complete answer to the question of when
the generalized McShane's identity holds for two generator
subgroups of $\SLTwoC$.

\vskip 10pt
\subsection{McShane-Bowditch identities for punctured torus
bundles}\label{ss:McShaneBowdithc} We next consider further
variations of the McShane-Bowditch identities. Recall that $\theta
\in \Gamma$ acts naturally on $\pi$, and hence on ${\mathcal X}$
where the action is given by
$$\theta([\rho])=[\rho \circ \theta^{-1}].$$

 Suppose
that $[\rho] \in {\mathcal X}$ is stabilized by an Anosov element
$\theta \in \Gamma^+$ (this corresponds to a hyperbolic element if
we identify $\Gamma^+$ with $\SLTwoZ$). So there exists $A \in
\SLTwoC$ such that for $\alpha \in \pi$,
$$\theta(\rho)(\alpha)=A\cdot \rho(\alpha)\cdot A^{-1}.$$
Note that  ${\rm tr}\,\rho(X)$ is well-defined on the classes $[X]
\in \mathscr C/\langle \theta \rangle $. Suppose further that
$[\rho]$ satisfies the \emph{relative} Bowditch Q-conditions on
${\mathscr C}/\langle \theta \rangle$, that is,
\begin{enumerate}
    \item ${\rm tr}\,\rho(X) \not \in [-2,2]$ for all $[X] \in \mathscr C/\langle \theta \rangle $;
    \item $|{\rm tr}\,\rho(X)| \le 2$ for only finitely many (possibly no)
    $[X] \in \mathscr C/\langle \theta \rangle $.
\end{enumerate}
Using the identification of $\Gamma^+$ with $\SLTwoZ$ and
${\mathscr C}$ with $\hat \Q \subset \hat \RR \cong {\mathscr
{PL}}$ in \S \ref{ss:simplecurves},  we get that the repelling and
attracting fixed points of $\theta$, $\mu_-, \mu_+ \in {\mathscr
{PL}}$ partitions ${\mathscr C}$ into two subsets ${\mathscr C}_L
\sqcup {\mathscr C}_R$ which are invariant under the action of
$\theta$. We have the following generalizations of the
McShane-Bowditch identities:
\begin{thm}\label{bundle} {\rm (}Theorems 5.6 and 5.9 of
\cite{tan-wong-zhang2004gMm}{\rm )} Suppose that $[\rho]$ is
stabilized by an Anosov element $\theta \in \Gamma^+$ and
satisfies the relative Bowditch Q-conditions as stated above. Then
\begin{eqnarray}\label{eqn:bundle1}
\sum_{[X] \in {\mathscr C}/\langle \theta
\rangle}\log\frac{e^{\nu}+e^{l(\rho(X))}}{e^{-\nu}+e^{l(\rho(X))}}
=0 \mod 2 \pi i,
\end{eqnarray}
and
\begin{eqnarray}\label{eqn:bundle2}
\sum_{[X] \in {\mathscr C}_L/\langle \theta
\rangle}\log\frac{e^{\nu}+e^{l(\rho(X))}}{e^{-\nu}+e^{l(\rho(X))}}
=\pm l(A) \mod 2 \pi i,
\end{eqnarray}
where the sums converge absolutely; and $l(A)$ is the complex
length of the conjugating element $A$ corresponding to $\theta$ as
described above, and the sign in (\ref{eqn:bundle2}) depends only
on our choice of orientations.
\end{thm}

\begin{rmk}
For type-preserving characters ($\kappa=-2$), the result is due to
Bowditch \cite{bowditch1997t}, where the summands of
(\ref{eqn:bundle1}) and (\ref{eqn:bundle2}) should be replaced
appropriately as in Remark \ref{rmk:1} by the summands of
McShane's original identity, and $l(A)$ in (\ref{eqn:bundle2})
should be replaced by $\lambda$, which has an interpretation as
the modulus of the cusp of an associated complete hyperbolic three
manifold. There are also similar identities in the case where
$\theta$ is reducible, that is, corresponds to a parabolic element
of $\SLTwoZ$, see \cite{tan-wong-zhang2004necsuf}.
\end{rmk}

The above result has applications to hyperbolic three manifolds.
Let $M$ be an orientable 3-manifold which fibers over the circle,
with the fiber a once-punctured torus, $\mathbb T$. Suppose that
the monodromy $\theta$ of $M$ is Anosov. By results of Thurston,
see \cite{thurston19??am} and \cite{thurston1978notes}, $M$ has a
complete finite-volume hyperbolic structure with a single cusp,
which can in turn be deformed to incomplete hyperbolic structures,
on which hyperbolic Dehn surgery can be performed to obtain
complete hyperbolic manifolds without cusps. Restricting the
holonomy representation to the fiber gives us characters which are
stabilized by $\theta$, and in the complete case, the relative
Bowdtich Q-conditions are satisfied, as for small deformations of
the complete structure to incomplete structures (see
\cite{tan-wong-zhang2004gMm}). The identities can be interpreted
as series identities for these (in)complete structures, involving
the complex lengths of certain geodesics corresponding to the
homotopy classes of essential simple closed curves on the fiber.
The quantity $\nu$ can be interpreted as half the complex length
of the meridian of the boundary torus, and $l(A)$ as the complex
length of a (suitably chosen) longitude of the boundary torus.

\vskip 5pt

\subsection{End Invariants}\label{ss:endinvariants}
We state some results concerning the end invariants of a character
$[\rho]$ in this subsection, these results can be found in
\cite{tan-wong-zhang2004endinvariants}.

\begin{defn}\label{def:endinvariant}(End invariants)

 An element $X \in \PL$ is  an end invariant
of $[\rho]$ if there exists $K>0$ and a sequence of distinct
elements $X_n \in {\mathscr C}$ such that $X_n \rightarrow X$ and
$|{\rm tr}\rho(X_n)|<K$ for all $n$.

\end{defn}
This definition generalizes the notion of a geometrically infinite
end for a discrete, faithful, type-preserving character. Denote by
$\EE$ the set of end invariants of $[\rho]$, this is a closed
subset of $\PL$ (see \cite{tan-wong-zhang2004endinvariants}).

\begin{thm}\label{thm:E=PL}
The set of end invariants $\EE$ is equal to $\PL$ if and only if
(i) $[\rho]$ is dihedral; or (ii) $[\rho]$ corresponds to a ${\rm
SU}(2)$ representation. Furthermore, if $\EE \neq \PL$, then $\EE$
has empty interior in $\PL$.

\end{thm}

\begin{thm}\label{thm:E=emptyset}
The set of end invariants $\EE$ is empty if and only if  $[\rho]$
satisfies the extended Bowditch Q-conditions.

\end{thm}
The next set of results  classify $\EE$ for real characters,
reducible characters ($\kappa=2$), imaginary characters and
discrete characters.
\begin{thm}\label{thm:EforReal}{\rm (}End invariants for real
characters{\rm )}. Suppose $[\rho] \in {\mathcal X}_{\kappa}$ is
real, with $\kappa \neq 2$. Then exactly one of the following must
hold:
\begin{enumerate}
\item [(a)] $\EE =\emptyset$, and  $\rho$ satisfies the extended
BQ-conditions.

\item[(b)] $\EE =\{\hat  X\}$ where $\hat X \in {\mathscr C}$,
$\rho$ is a $\SLTwoR$ representation, ${\rm tr}\,\rho(\hat X) \in
(-2,2)$, and ${\rm tr}\,\rho( X) \not\in (-2,2)$ for all $X \in
{\mathscr C} \setminus \{ \hat X \}$.

\item[(c)] $\EE$ is a Cantor subset of $\PL$,  $\rho$ is a
$\SLTwoR$ representation,   ${\rm tr}\,\rho(X) \in (-2,2)$ for at
least two distinct $X \in {\mathscr C}$, and ${\rm tr}\,\rho(Y)
\not \in (-2,2) \cup \{\pm \sqrt{\kappa+2}\}$ for some element $Y
\in {\mathscr C}$.

\item[(d)] $\EE=\PL$, and $\rho$ satisfies the conditions of
Theorem \ref{thm:E=PL}, that is, $\rho$ is the dihedral
representation or a $\SUTwo$ representation.
\end{enumerate}

Furthermore, case (a) occurs only when $ \kappa \in (-\infty,
2)\cup [\,18, \infty)$; case (b) when $ \kappa \in [\,6, \infty)$;
case (c) when $ \kappa \in (2,\infty) $; and case (d) when $
\kappa \in [-2,2)\cup (2,\infty]$.
\end{thm}

\begin{thm}\label{thm:Eforreducible}{\rm (}End Invariants for reducible
characters{\rm )}.\\
 For $ [\rho] \in {\mathcal X}_2$,   $\EE=\{X_0\}$ or
$\PL$.
 Furthermore, in the first case, if $X_0 \in
{\mathscr C}$, then ${\rm tr}\, \rho(X_0) \in [-2,2]$ and ${\rm
tr}\, \rho(X) \not\in [-2,2]$ for all $X \in {\mathscr C}\setminus
\{X_0\}$, if $X_0 \not\in {\mathscr C}$, then ${\rm tr}\, \rho(X)
\not\in [-2,2]$ for all $X \in {\mathscr C}$; and in the second
case, ${\rm tr}\, \rho(X) \in [-2,2]$ for all $X \in {\mathscr
C}$.

\end{thm}

The following theorem gives a partial classification for imaginary
characters.

\begin{thm}\label{thm:EforImaginary}{\rm (}End Invariants for imaginary
characters{\rm )}. Suppose that $[\rho]$ is imaginary.
\begin{enumerate}
\item [(i)] $\kappa =-2$: For $ [\rho] \in {\mathcal X}_{-2}$,
$\EE$ is either a Cantor subset of $\PL$, or consists of a single
element $X$ in ${\mathscr C}$. In the latter case, ${\rm
tr}\,\rho(X)=0$ and $[\rho]$ is equivalent under the action of
$\Gamma$ to a character corresponding to the triple $(0, x, ix)$
where $x \in \RR$ satisfies $|x| \ge 2$.

\item[(ii)] $-14 \le\kappa < 2$: For $ [\rho] \in {\mathcal
X}_{\kappa}$, $\EE$ is either a Cantor subset of $\PL$, or
consists of a single element $X$ in ${\mathscr C}$.

\item[(iii)] $\kappa <-14$: For $ [\rho] \in {\mathcal
X}_{\kappa}$, $\EE$ is  a Cantor subset of $\PL$; consists of a
single element $X$ in ${\mathscr C}$; or is empty.

\end{enumerate}

\end{thm}

Finally, we say that a character is discrete if the set $\{ {\rm
tr}\, \rho(X) ~|~ X \in {\mathscr C}\}$ is discrete in $\CC$. We
have the following result.

\begin{thm}\label{thm:Efordiscrete}
For a discrete $ [\rho] \in {\mathcal X}$, if  $\EE$ has at least
three elements, then $\EE$ is either a Cantor subset of $\PL$ or
all of $\PL$.

\end{thm}


\vskip 20pt
\section{{\bf Further directions and related results}}\label{s:further directions}
\vskip 20pt

A natural question to pose is whether the results above extend to
general surfaces, and especially to closed surfaces without
boundary (where a suitable version of McShane's identity is still
lacking). The identities (\ref{eqn:TWZ}) and (\ref{eqn:McShane})
have been  been generalized by McShane himself to hyperbolic
surfaces with cusps \cite{mcshane1998im},to hyperbolic surfaces
with cusps and/or geodesic boundary components by Mirzakhani
\cite{mirzakhani2004preprint}, to hyperbolic surfaces with cusps,
geodesic boundary and/or conical singularities, as well as to
classical Schottky groups by the authors in
\cite{tan-wong-zhang2004cone-surfaces},
\cite{tan-wong-zhang2004schottky}. Further refinements and
analogous results for punctured surface bundles over the circle
have also been obtained by Akiyoshi, Miyachi and Sakuma in
\cite{akiyoshi-miyachi-sakuma2004cm355}
\cite{akiyoshi-miyachi-sakuma2004preprint}. The methods used in
the above cited works follow closely that used by McShane in his
original proof and differs markedly from those used to prove the
results announced in this note, which are modelled on Bowditch's
proof. In particular, a general set of conditions, equivalent to
the Bowditch Q-conditions, for the identities to hold  is still
lacking. The Bowditch method does extend fairly naturally to the
case of the four holed sphere, see \cite{goldman-tan-zhang2005}.
The one-holed torus and the four-holed sphere are the natural
building blocks for more complicated surfaces, according to
Grothendieck's reconstruction principle \cite{grothendieck1997},
see also Luo's work on ${\rm SL}(2)$ character varieties of
surfaces in \cite{luo1999}. It would be interesting to see if this
method can be used to give an independent proof of the identities
for characters of general surfaces, and to provide necessary and
sufficient conditions for the identities to hold, and also to shed
some light on the dynamics of the mapping class group action on
the $\SLTwoC$ character varieties for these general surfaces.

Mirzakhani found some beautiful applications of McShane's
identity, in particular, she used it to compute the Weil-Petersson
volume of the moduli space of bordered surfaces in
\cite{mirzakhani2004preprint}, see also
\cite{mirzakhani_2_2004preprint} and
\cite{mirzakhani2004_3_preprint} for other striking applications.
Theorem \ref{thm:TWZ} gives a natural definition of a moduli space
for (relative) $\SLTwoC$ characters and it would be interesting to
know if the volume can be computed for this moduli space, with
respect to the Poisson (complex-symplectic) structure which is
invariant under the mapping class group action.

\vskip 5pt There are also interesting generalizations concerning
the end invariants of $[\rho]$. The results stated in \S
\ref{ss:endinvariants} can be considered as evidence towards the
following conjecture, which is a refinement and generalization of
the suggestion by Bowditch in \cite{bowditch1998plms} that for a
generic $[\rho] \in {\mathcal X}_{-2}$ not satisfying the
BQ-conditions, $\EE$ should be a Cantor set.

\begin{conj}\label{conj:cantorset}
Suppose that $\EE$ has more than two elements. Then either
$\EE=\PL$ or $\EE$ is a Cantor subset of $\PL$.
\end{conj}

\vskip 5pt

There is also a generalization of the ending lamination conjecture
 for $\SLTwoC$ characters  (as Bowditch conjectured for the
$\kappa=-2$ case) which can be stated as follows:

\begin{conj}\label{conj:endlamination}
Suppose that $[\rho], [\rho'] \in {\mathcal X}_{\kappa}$ are such
that $\EE={\mathcal E}(\rho')$, $\EE$ has at least two elements,
and $\EE \neq\PL$. Then $[\rho]=[\rho']$.

\end{conj}

Finally, we say a few words about the  generalizations to
arbitrary surfaces. The definition of $\EE$ can be extended
without much difficulty. The case of the four-holed sphere is
similar and the techniques given here should give similar results
in that case. In other cases, $\PL$ is homeomorphic to the sphere
$S^n$ for some $n \ge 2$ and a possible generalization of Theorem
\ref{thm:E=PL} is that $\EE$ has either full measure or measure
zero. A possible generalization of Conjecture \ref{conj:cantorset}
would be that $\EE$ is perfect, if it contains more than two
elements. However, we do not have any insights into  these more
general cases.

\end{document}